\documentclass[a4paper,11pt,twoside]{article}
\usepackage{amsmath,amssymb,amsthm}

\newcommand{\defeq}{:=}
\renewcommand{\k}{\Bbbk}
\newcommand{\tens}{\otimes}
\newcommand{\xd}{\mathrm{d}}
\newcommand{\lalg}[1]{\mathfrak{#1}}
\newcommand{\ealg}{\mathsf{U}}
\newcommand{\func}{\mathcal{C}}
\DeclareMathOperator{\id}{id}

\DeclareMathOperator{\cou}{\epsilon}
\DeclareMathOperator{\cop}{\Delta}
\DeclareMathOperator{\antip}{\mathrm{S}}
\newcommand{\act}{\triangleright}
\newcommand{\ract}{\triangleleft}
\renewcommand{\i}[1]{{}_{\scriptscriptstyle(#1)}}
\newcommand{\iu}[1]{{}_{\scriptscriptstyle(\underline #1)}}
\newcommand{\catmod}[4]{{}^{#1}_{#2}{\mathcal{M}}^{#3}_{#4}}
\newcommand{\catmodx}[4]{{}^{#1}_{#2}{\dot{\mathcal{M}}}^{#3}_{#4}}
\newcommand{\ftr}[1]{\mathcal{#1}}

\theoremstyle{plain}
\newtheorem{prop}{Proposition}[section]

\newtheorem{thm}[prop]{Theorem}

\newtheorem{dfn}[prop]{Definition}

\title{\textbf{Structure Theorem for Covariant Bundles on Quantum
Homogeneous Spaces}}
\author{Robert Oeckl\footnote{email: oeckl@cpt.univ-mrs.fr}\\ \\
Centre de Physique Th\'eorique, CNRS Luminy,\\ 13288 Marseille
cedex 9, France}
\date{CPT-2001/P.4240\\ 6 August 2001}

\begin{document}
\maketitle

\begin{abstract}
The natural generalization of the notion of bundle in
quantum geometry is that of bimodule. If the base space
has quantum group symmetries one is particularly interested
in bimodules covariant (equivariant) under these symmetries.
Most attention has so far been focused on the case with maximal
symmetry --
where the base space is a quantum group and the bimodules are
bicovariant. 
The structure of bicovariant bimodules is well understood through
their correspondence with crossed modules.

We investigate the ``next best'' case -- where the base space is a
quantum homogeneous space and the bimodules are covariant.
We present a structure theorem that resembles the one for
bicovariant bimodules.
Thus, there is a correspondence between covariant
bimodules and a new kind of ``crossed'' modules which we define.
The latter are attached to the pair of quantum groups which defines the
quantum homogeneous space.

We apply our structure theorem to differential
calculi on quantum homogeneous spaces 
and discuss a related notion of induced differential calculus.
\end{abstract}

\section{Preliminaries}

We start by introducing notation and reviewing some relevant
definitions. Thus, coproduct, counit and antipode of a Hopf algebra
are denoted $\cop,\cou,\antip$ respectively. We use Sweedler's
notation (with implicit summation) $\cop h=h\i1\tens h\i2$ for the
coproduct. A similar notation serves for left coactions $v\mapsto
v\i1\tens v\iu2$, and correspondingly for right ones
(the component remaining in the comodule is underlined).
Throughout we work over a field $\k$.

Let $H$ be a Hopf algebra. We denote the category of left $H$-modules
by $\catmod{}{H}{}{}$, the category of left $H$-comodules by
$\catmod{H}{}{}{}$, and correspondingly for the right hand side
versions. Furthermore, for modules that carry several (co)actions
which mutually commute we use the obvious notation for the
category. (Such modules are also called \emph{Hopf modules}.)
E.g., for left $H$-module right $H$-comodules such that both
structures are compatible we would write $\catmod{}{H}{H}{}$.
A module with a compatible comodule structure is also called a
\emph{covariant} module. Compatible left and right (co)module
structures are called \emph{bi(co)module}.

We consider a second type of module which is called \emph{crossed
module} (or \emph{Yetter-Drinfeld module}). Let $H$ be a Hopf
algebra. A right crossed $H$-module $V$ is a right $H$-module and
right $H$-comodule such that the following condition holds:
\[
 v\iu1 \ract h\i1\tens v\i2 h\i2= 
 (v\ract h\i2)\iu1\tens h\i1 (v\ract h\i2)\i2 \quad \forall h\in H,
 v\in V
\]
We denote the category of such modules by $\catmodx{}{}{H}{H}$.
There is also a corresponding left handed version.

The structure theorem for bicovariant bimodules (playing the role of
bicovariant bundles over a quantum group) can be formulated as
follows (in its right handed version). This result is implicit to some
extent in \cite{Wor:calculi}. A  complete formulation was given
in \cite{Sch:modules}.

\begin{thm}
\label{thm:bicov}
Let $H$ be a Hopf algebra. The categories $\catmod{H}{H}{H}{H}$ and
$\catmodx{}{}{H}{H}$ are equivalent.
\end{thm}

The equivalence is given in one direction by the functor
$\catmod{H}{H}{H}{H}\to \catmodx{}{}{H}{H}$ defined by
$E\mapsto {}^H E\defeq \{e\in E : e\i1\tens e\iu2 = 1\tens e\}$. ${}^H
E$ inherits the right $H$-comodule structure from $E$ and is equipped
with the new right $H$-module structure $e \tilde{\ract} h \defeq \antip h\i1
\act e \ract h\i2$. Conversely, the inverse functor
$\catmodx{}{}{H}{H}\to \catmod{H}{H}{H}{H}$ is given by $X\mapsto
H\tens X$. Here, the left module and comodule structure of $H\tens X$
are the regular ones of $H$ while the right structures are the tensor
product ones.

We will be interested in quantum homogeneous spaces defined as
follows.
\begin{dfn}
Let $\pi:P\to H$ be a surjection of Hopf algebras.
Then the left $P$-comodule algebra $B\defeq P^H = \{p\in P : p\i1\tens
\pi(p\i2) = p\tens 1\}$ is called a \emph{quantum homogeneous space}.

The triple $(P,B,H)$ is said to satisfy the \emph{Hopf-Galois} property
if the map $\chi:P\tens_B P\to P\tens H$ given by
$\chi=(\cdot\tens\pi)\circ(\id\tens\cop)$ is injective (in addition
to being surjective).
\end{dfn}

Note that the Hopf-Galois
condition is automatically
satisfied if $H$ is cosemisimple (and thus also has invertible
antipode).
This follows from \cite{Sch:principalhom}. (Apply Remark 3.3.(2) to
the integral and use Remark 3.3.(1) in Theorem 3.5.)

A bundle structure of prime importance in differential geometry is
the (co)tangent bundle. A noncommutative generalization of this
notion (together with the exterior derivative of functions)
is captured by the notion of differential calculus given as follows.
\begin{dfn}
Let $B$ be an algebra. A \emph{differential calculus} $\Omega$ over
$B$ is a
$B$-bimodule with a linear map $\xd:B\to\Omega$ such that (a) the
Leibniz rule $\xd(ab)=a\xd(b)+\xd(a) b$ is satisfied and (b) the map
$B\tens B\to\Omega: a\tens b\mapsto a\xd b$ is surjective.
\end{dfn}

A basic result about differential calculi is the following (see e.g.\
\cite{Wor:calculi}).
\begin{prop}
\label{prop:univ}
Let $B$ be an algebra.
The \emph{universal differential calculus} over $B$ is given by
$\tilde{B}\defeq\ker\cdot\subset B\tens B$ with left and right $B$-module
structures  given by multiplication of the left respectively right
component. The exterior derivative $\xd: B\to \tilde{B}$ is given by
$\xd b=1\tens b-b\tens 1$.
Any differential calculus over $B$ can be identified with a quotient
of $\tilde{B}$ by a subbimodule.
\end{prop}

If the base space $B$ has extra symmetries as in the case of a quantum
homogeneous space or even a quantum group it is natural to demand these
symmetries also from the differential calculus (as in the commutative
situation). This leads to the
obvious notions of covariant or bicovariant differential calculus.
Proposition~\ref{prop:univ} remains valid if quotient is understood to
mean quotient by a subbimodule which is (bi)covariant.

\section{Induced Differential Calculi -- Motivation}

In this section we construct a differential calculus on a quantum
homogeneous space from a given one on the symmetry quantum group.
In fact, this is nothing but the construction of the
cotangent bundle on a homogeneous space from the cotangent bundle on
the symmetry group -- but formulated in a way that generalizes to the
noncommutative case. We use a notation that is intended to remind
the reader of the differential geometric origin.

Recall that (by application of the Structure Theorem~\ref{thm:bicov})
a bicovariant differential calculus over a quantum group
$\func(G)$ is given by a bicovariant bimodule
$\Gamma(T^* G)=\func(G)\tens T^*_e G$
(classically the space of
sections of the cotangent bundle over the Lie group $G$)
\cite{Wor:calculi}. $T^*_e G$ is the right crossed
$\func(G)$-module of left-invariant 1-forms, which corresponds
classically to the cotangent space at the identity of $G$.
$T^*_e G$ is a quotient $\func(G)^+/I$ of
$\func(G)^{+}\defeq\ker\cou\subset\func(G)$ as a right crossed
$\func(G)$-module via the right regular action and right adjoint
coaction. The exterior derivative $\xd:\func(G)\to \func(G)\tens T^*_e
G$ is determined by $f\mapsto f\i1\tens f\i2 - f\tens 1$.
In the classical case $I=(\ker\cou)^2$. Then $I$
is the annihilator of $\lalg{g}\subset \ealg(\lalg{g})$ (with
$\lalg{g}$ the Lie algebra of $G$)
in the pairing
of $\func(G)$ with $\ealg(\lalg{g})$ and thus
$\func(G)^+/I\cong \lalg{g}^*$.

\begin{prop}[Induced differentials on homogeneous spaces]
Let $\pi:\func(G)\to \func(H)$ be a surjection of Hopf algebras
with $\func(M)\defeq \func(G)^{\func(H)}$.
Let $\Gamma(T^* G)=\func(G)\tens T^*_e G$ be a bicovariant
differential calculus on $\func(G)$.
We obtain a corresponding differential calculus on the homogeneous
space $\func(M)$ in two steps.

First, we restrict the
cotangent space at each point to those forms that are annihilated
by the vector fields generated by the right translations of $H$.
Thus, we define
\[
 T^*_e M\defeq \func(M)\cap (\func(G)^+/I)
= \func(M)^+/(\func(M)^+\cap I)
\]
with $\func(M)^+\defeq \ker\cou\subset\func(M)$.
While $T^*_e M$ does not carry a right $\func(G)$ coaction anymore
it does inherit from $T^*_e G$ the induced right coaction of $\func(H)$.
Furthermore, it carries a right $\func(M)$ action, the restriction of
the right $\func(G)$ action on $T^*_e G$.

Now, the second step consists in restricting the so formed ``bundle
over $G$'' to a ``bundle over $M$''. This is accomplished
by going to the $\func(H)$-invariant subspace
\[
 \Gamma(T^* M)\defeq (\func(G)\tens T^*_e M)^{\func(H)} .
\]
This is now a left $\func(G)$-covariant $\func(M)$-bimodule.
$d:\func(G)\to \Gamma(T^* G)$ descends to a map
$d:\func(M)\to \Gamma(T^* M)$.
The classical case recovers the usual differential calculus on $M$.
\end{prop}
\begin{proof}
The induced right adjoint coaction of $\func(H)$ on $\func(G)$ is
closed on the subspace $\func(M)$:
\[
 a\mapsto a\i2\tens \pi(\antip a\i1 a\i3)
 = a\i2 \tens \pi(\antip a\i1) \in \func(M)\tens \func(H)
\]
for $a\in \func(M)$.

That $d$ descends follows for step 1 from
$\cop \func(M)\subseteq \func(G)\tens\func(M)$ and for step 2 from
the right $\func(H)$-invariance of $\cop \func(M)$:
\[
 a\i1\tens a\i2\mapsto a\i1\tens a\i4\tens \pi(a\i2)
 \pi(\antip a\i3 a\i5) = a\i1\tens a\i2\tens 1
\]
for $a\in \func(M)$.
\end{proof}

\section{Structure Theorem}

The structure found for induced differential calculi in the
previous section naturally leads to the question whether this
structure is generic.
We present here our main result, a
structure theorem for covariant bimodules (i.e., covariant bundles)
over quantum homogeneous spaces in analogy to Theorem~\ref{thm:bicov}.
This answers the question in the
affirmative. Finally, we reapply the structure theorem to differential
calculi.

\begin{dfn}
Let $P\to H$ be a surjection of Hopf algebras, $B\defeq P^H$.
Then a right $B$-module and right $H$-comodule $X$ is called
\emph{crossed} iff
\[
 x\iu1 \ract b \tens x\i2 = (x\ract b\i2)\iu1 \tens 
 \pi(b\i1) (x\ract b\i2)\i2
\]
for all $x\in X, b\in B$. We denote the category of such objects by
$\catmodx{}{}{H}{B}$.
\end{dfn}

Note that if $H$ has invertible antipode the right $H$-coaction can be
converted to a left $H$-coaction $x\mapsto \antip^{-1} x\i2\tens
x\iu1$ and $\catmodx{}{}{H}{B}\cong\catmod{H}{}{}{B}$.

\begin{thm}
\label{thm:cormod1}
Let $P\to H$ be a surjection of Hopf algebras, $B\defeq P^H$.
Then, the categories $\catmodx{}{}{H}{B}$ and $\catmod{P}{P}{H}{B}$
are equivalent.
\end{thm}
\begin{proof}
For $X\in \catmodx{}{}{H}{B}$ consider the tensor product $P\tens
X$. Equip it with the left (co)module structures of $P$ and
the right (co)module structures of the tensor product.
One checks that this makes $P\tens X$ an object in $\catmod{P}{P}{H}{B}$.
For a morphism $f:X\to X'$ in $\catmodx{}{}{H}{B}$ the map $\id\tens
f:P\tens X\to P\tens X'$ is a morphism in $\catmod{P}{P}{H}{B}$. This
defines a functor $\catmodx{}{}{H}{B}\to \catmod{P}{P}{H}{B}$.

For $Y\in\catmod{P}{P}{H}{B}$ consider its left
$P$-invariant subspace ${}^P Y$. We equip it with a new right action of $B$
via $y\,\tilde{\ract}\, b\defeq \antip b\i1 \act y \ract b\i2$. This makes
it with the inherited right coaction of $H$ an object in
$\catmodx{}{}{H}{B}$. A morphism $g:Y\to Y'$ in
$\catmod{P}{P}{H}{B}$ induces a morphism $\tilde{g}: {}^P Y\to {}^P
Y'$ in $\catmodx{}{}{H}{B}$ by restriction.
This defines a functor $\catmod{P}{P}{H}{B}\to\catmodx{}{}{H}{B}$.

Finally, we check that the two functors are mutually inverse.
While clearly ${}^P (P\tens X)\cong X$ the isomorphism
$Y\cong P\tens {}^P Y$ is given by
$y\mapsto y\i1\tens \antip y\i2\act y\iu3$
with inverse $p\tens y\mapsto p\act y$. To check the inverseness for
morphisms is straightforward and left to the reader.
\end{proof}

\begin{thm}
\label{thm:cormod2}
Let $P\to H$ be a surjection of Hopf algebras, $B\defeq P^H$.
Then, there are functors $\ftr{F}:\catmod{P}{B}{}{B}\to
\catmod{P}{P}{H}{B}$ and $\ftr{G}:\catmod{P}{P}{H}{B}\to
\catmod{P}{B}{}{B}$ such that $\ftr{G}\circ\ftr{F}$ is the identity.

If furthermore $H$ has invertible antipode and $(P,B,H)$ is
Hopf-Galois, then also
$\ftr{F}\circ\ftr{G}$ is the identity and the categories are thus
equivalent.
\end{thm}
\begin{proof}
For $E\in\catmod{P}{B}{}{B}$ consider the tensor product
$P\tens E$. We equip it with the left $P$-comodule
structure as a tensor product, the left
$P$-module and right $H$-comodule structure of $P$ and
the right $B$-module structure of $E$. These structures descend to the
quotient $P\tens_B E$ and make it an object in
$\catmod{P}{P}{H}{B}$. A morphism $h:E\to E'$ in $\catmod{P}{B}{}{B}$
defines a map $\id\tens h:P\tens E\to P\tens E'$ which induces a
morphism $\tilde{h}:P\tens_B E\to P\tens_B E'$ in
$\catmod{P}{P}{H}{B}$. This defines the functor
$\ftr{F}:\catmod{P}{B}{}{B}\to \catmod{P}{P}{H}{B}$.

Given $Y\in\catmod{P}{P}{H}{B}$ consider the right
$H$-invariant subspace $Y^H$. The left $P$-comodule and right
$B$-module structures descend while the left $P$-module structure
only survives as a left $B$-module structure. Thus,
$Y^H\in\catmod{P}{B}{}{B}$. A morphism $g:Y\to Y'$ in
$\catmod{P}{P}{H}{B}$ clearly gives rise to a morphism
$\tilde{g}:Y^H\to {Y'}^H$ in $\catmod{P}{B}{}{B}$
by restriction. This defines the functor $\ftr{G}:\catmod{P}{P}{H}{B}\to
\catmod{P}{B}{}{B}$.

Next, we check that $\ftr{G}\circ\ftr{F}=\id$.
Starting with $E$ we obtain the corresponding object
$\ftr{F}(E)=P\tens_B E$ and $(\ftr{G}\circ\ftr{F})(E)=(P\tens_B
E)^H$. But this is $B\tens_B E$ (as taking invariant subspace and
quotient of the tensor product commute by construction) which in turn
is canonically isomorphic to $E$. To check $\ftr{G}\circ\ftr{F}=\id$
on morphisms is straightforward and left to the reader.

We assume now further that $H$ has invertible antipode and that
$(P,B,H)$ is Hopf-Galois.
For $X\in\catmodx{}{}{H}{B}$
consider the map $\chi\tens\id: P\tens_B P\tens X\to P\tens H\tens
X$. We define a right coaction of $H$ on both
sides as the tensor product one on $P\tens X$ and $H\tens X$
respectively. (This definition behaves well with respect to the tensor
product $\tens_B$.)
It commutes with $\chi\tens\id$ which thus restricts
to a map on the invariant subspaces under this coaction
$\tilde{\chi}:P\tens_B (P\tens X)^H\to P\tens (H\tens X)^H$.
As $(P,B,H)$ is Hopf-Galois, $\chi$ is a bijection and so are
$\chi\tens\id$ and $\tilde{\chi}$.

Now the map $(H\tens X)^H\to X$ given by $\cou\tens\id$ is a bijection
since the antipode of $H$ is invertible. Its inverse is given by
$x\mapsto \antip^{-1} x\i2\tens x\iu1$. Thus, we obtain bijections
$P\tens_B (P\tens X)^H\to P\tens (H\tens X)^H\to P\tens X$. Using
Theorem~\ref{thm:cormod1} this gives rise to a bijection
$(\ftr{F}\circ\ftr{G})(Y)=P\tens_B Y^H\to Y$ for any
$Y\in\catmod{P}{P}{H}{B}$. One easily checks that this is an
isomorphism with respect to the relevant (co)module structures.
Thus, $\ftr{F}\circ\ftr{G}$ is the identity on objects. To check that
it is the identity on maps is now straightforward and left to the reader.
\end{proof}

\begin{prop}
\label{prop:dchom}
Let $\pi:P\to H$ be a surjection of Hopf algebras, $B\defeq P^H$.
Let $B^+\defeq \ker\cou\subset B$.
Then, each left $P$-covariant differential calculus on $B$ corresponds
to a crossed submodule $I\subset B^+\in\catmodx{}{}{H}{B}$ via 
the right regular action and the coaction
$b\mapsto b\i2 \tens \pi(\antip b\i1)$.

If furthermore $H$ has invertible antipode and $(P,B,H)$ is
Hopf-Galois, then the correspondence is one-to-one.
\end{prop}
\begin{proof}
We use the fact that any differential calculus is a quotient
of the universal one (Proposition~\ref{prop:univ}) and apply the
correspondences of the previous
theorems. On $B$ the universal calculus is given by the subspace
$\ker\cdot\subset B\tens B$ with $\xd:B\to B\tens B$ defined
by $b\mapsto 1\tens b - b\tens 1$. For simplicity we start by
considering the whole space $B\tens B$.
It is a left $P$-covariant $B$-bimodule (i.e.\ an object in
$\catmod{P}{B}{}{B}$) by the tensor product
coaction and the left and right regular actions of $B$ on the
left and right component respectively. With Theorem~\ref{thm:cormod2}
we obtain $P\tens_B (B\tens B)\cong P\tens B$ as
an object in $\catmod{P}{P}{H}{B}$.
This in turn corresponds to $^P (P\tens B)$ as an object in
$\catmodx{}{}{H}{B}$ according
to Theorem~\ref{thm:cormod1}. This in turn we can identify with $B$
via the map $p\tens b\mapsto \cou(p) b$ and its inverse
$b\mapsto \antip b\i1\tens b\i2$.
The induced module structure on $B$ is the right regular action while
the right $H$-comodule structure is $b\mapsto b\i2\tens\pi(\antip
b\i1)$.
By applying the inverse functors of Theorems~\ref{thm:cormod2} and
\ref{thm:cormod1} to this $B$ we obtain an isomorphism
$B\tens B\to (P\tens B)^H$ in $\catmod{P}{B}{}{B}$
given by $b\tens c\mapsto b c\i1\tens
c\i2$ with inverse $p\tens b\mapsto p \antip b\i1\tens b\i2$.

Now the subspace $\ker\cdot\subset B\tens B$ on the left hand side
corresponds to $(P\tens B^+)^H$ on the right hand side.
The differential map $d:B\to (P\tens B^+)^H$ is
$b\mapsto b\i1\tens b\i2 - b\tens 1$.

Since left $P$-covariant differential calculi on $B$ correspond to
quotients of $\ker\cdot\subset B\tens B$ in $\catmod{P}{B}{}{B}$, by
the above correspondence these in turn correspond to
quotients (and thus crossed submodules $I$) of $B^+$ in
$\catmodx{}{}{H}{B}$.

In general each differential calculus corresponds to a certain such
crossed submodule (as the composition $\ftr{G}\circ\ftr{F}=\id$ in
Theorem~\ref{thm:cormod2}).
If the additional condition of invertibility of the
antipode of $H$ and the Hopf-Galois for $(P,B,H)$ is satisfied the
converse is also true, giving rise to a one-to-one correspondence
(as also $\ftr{F}\circ\ftr{G}=\id$ in Theorem~\ref{thm:cormod2}).
\end{proof}

Note that a result similar to Proposition~\ref{prop:dchom} for the
more general case where $B$ is a coideal subalgebra of $P$ was found
recently by direct means \cite{Her:derqg}.

\section*{Acknowledgements}

I would like to thank Istvan Heckenberger for interesting
discussions and useful comments on the manuscript.
This work was supported by a NATO fellowship grant.

\end{document}